\catcode`\@=11

\font\tenmsa=msam10 \font\sevenmsa=msam7 \font\fivemsa=msam5
\newfam\msafam
\textfont\msafam=\tenmsa \scriptfont\msafam=\sevenmsa
\scriptscriptfont\msafam=\fivemsa
\def\hexnumber@#1{\ifnum#1<10 \number#1\else \ifnum#1=10
A\else\ifnum#1=11
 B\else\ifnum#1=12 C\else \ifnum#1=13 D\else\ifnum#1=14
E\else\ifnum#1=15
 F\fi\fi\fi\fi\fi\fi\fi}
\def\msa@{\hexnumber@\msafam}

\mathchardef\blacktriangleright="3\msa@49
\mathchardef\blacktriangleleft="3\msa@4A

\catcode`\@=\active

\def\CC{\hbox{{$\cal C$}}}
\def\CV{\hbox{{$\cal V$}}}
\def\CM{\hbox{{$\cal M$}}}
\def\CL{\hbox{{$\cal L$}}}
\def\CR{\hbox{{$\cal R$}}}

\def\Nat{{\rm Nat}}
\def\Aut{{\rm Aut}}
\def\eps{{\epsilon}}
\def\lcross{{>\!\!\!\triangleleft}}
\def\rcocross{{\blacktriangleright\!\!<}}
\def\dcross{{\bowtie}}
\def\rbiprod{{\cdot\kern-.33em\triangleright\!\!\!<}}
\def\lbiprod{{>\!\!\!\triangleleft\kern-.33em\cdot}}

\def\cora{\blacktriangleleft}
\def\tens{\mathop{\otimes}}
\def\la{{\triangleright}}
\def\ra{{\triangleleft}}
\def\isom{{\cong}}
\def\Mor{{\rm Mor}}
\def\ev{{\rm ev}}
\def\coev{{\rm coev}}
\def\id{{\rm id}}
\def\<{\langle}
\def\>{\rangle}
\def\vecv{{\bf v}}
\def\vecu{{\bf u}}
\def\und#1{{\underline {#1}}}
\def\note#1{}
\def\proof{\goodbreak\noindent{\bf Proof\quad}}
\def\lform{\hbox{$\sqcup$}\llap{\hbox{$\sqcap$}}}
\def\endproof{{\ $\lform$}\bigskip }

\documentstyle[epsf,11pt]{article}  

\newtheorem{lemma}{Lemma}[section]
\newtheorem{propos}[lemma]{Proposition}
\newtheorem{theorem}[lemma]{Theorem}

\begin{document}\baselineskip 18pt

{\ }\hskip 4.7in   Damtp/97-140
\vspace{.2in}

\begin{center} {\Large SOME REMARKS ON BRAIDED GROUP RECONSTRUCTION}\\
{\Large AND BRAIDED DOUBLES}\\ 
\baselineskip 13pt{\ }
{\ }\\ Shahn Majid\footnote{Royal Society University Research
Fellow and Fellow of Pembroke College, Cambridge, England.}\\ {\
}\\ Department of Applied Mathematics and Theoretical Physics\\
University of Cambridge, Cambridge CB3 9EW, UK\\
www.damtp.cam.ac.uk/user/majid
\end{center}
\begin{center}
Revised August, 1998
\end{center}
\vspace{10pt}

\begin{quote}\baselineskip 13pt {\bf ABSTRACT} The cross coproduct 
braided group $\Aut(\CC)\rcocross B$ is obtained by Tannaka-Krein 
reconstruction from $\CC^B\to \CC$ for a braided group $B$ in braided 
category $\CC$. We apply this
construction to obtain partial solutions to two problems in braided
group theory, namely the tensor problem and the braided double. We
obtain $\Aut(\CC)\rcocross\Aut(\CC)$ $\isom \Aut(\CC)\lcross\Aut(\CC)$
and higher braided group `spin chains'. The example of the braided
group $B(R)\lcross
B(R)\lcross\cdots\lcross B(R)$ is described explicitly by R-matrix
relations. We also obtain $\Aut(\CC)\rcocross \Aut(\CC)^*$ as a
dual quasitriangular `codouble' braided group by reconstruction
from the dual category $\CC^\circ\to\CC$. General braided double
crossproducts $B\dcross C$ are also considered.
\end{quote}\baselineskip 18pt

\section{Introduction}

Recently there has been a lot of interest in the braided version of
the Tannaka-Krein reconstruction theorem proven by the author in
\cite{Ma:bg}. Here one considers a monoidal functor $F:\CC\to \CV$
  and reconstructs (under certain
representability assumptions) a braided group or Hopf algebra in
the braided category $\CV$, denoted $\Aut(\CC,F,\CV)$. This is how
braided groups were first introduced (in \cite{Ma:bg}). By now
there is a rich and extensive theory of braided groups, see for
example \cite{Ma:bra}\cite{Ma:tra}\cite{Ma:bos}\cite{Ma:rep}
\cite{Ma:cat}\cite{Ma:exa}\cite{Ma:lin}\cite{Ma:skl}\cite{Ma:lie}
\cite{Ma:diag}\cite{Ma:dbos},
and \cite{Ma:introm}\cite{Ma:book} for reviews.

In this paper we consider some examples of braided
reconstruction. Let $B$ be a braided group in a category $\CC$ and
$\CC^B$ the category of braided comodules of $B$. In Section~1 we
study the solution of the reconstruction problem for the forgetful
functor
\[ \CC^B\to \CC,\]
namely the prebosonisation braided group cross coproduct
$\Aut(\CC)\rcocross B$ introduced in \cite{Ma:bos} (there in the
module version rather than comodule version but the reversal of
arrows is routine). It is the prebosonisation of $B$ because in the
case $\CC=\CM^H$ ($H$ dual quasitriangular) it is related by
transmutation to the bosonisation $H\rbiprod B$. All this was known
since 1991 from the bosonisation theory
\footnote{The coalgebra part of this reconstruction problem
with $\CC=\CM^H$ was recently considered in the preprint 
of \cite{Par:rec}. The full solution (the entire braided group structure) 
and the fact that it was already known from \cite{Ma:bos} was 
pointed out by the author in a letter to Pareigis during his 
preparation of the final version of \cite{Par:rec}.}; to this
we add now the observation that when $B=\Aut(\CC)$ itself,
\[ \Aut(\CC)\rcocross \Aut(\CC)\isom \Aut(\CC)\lcross\Aut(\CC),\]
i.e. in some sense should be viewed as the closest one may come to
something in between like $\Aut(\CC)\tens\Aut(\CC)$, having an
equal description either as a tensor product algebra with cross
coalgebra or a tensor product coalgebra with a cross algebra. (We
recall that there is no general tensor product of braided groups in
a given braided category; the tensor product algebra and coalgebra
themselves do not fit together due to `tangling up' when the
category is truly braided.) The construction can be iterated and
leads to concrete R-matrix formulae for a $n$-fold `spin chain'
braided group $B(R)\lcross B(R)\lcross\cdots\lcross B(R)$ where
$B(R)$ are the braided matrices \cite{Ma:exa}.

In Section~3 we solve the reconstruction problem for the forgetful
functor
\[ \CC^\circ\to \CC\]
where $\CC^\circ$ is the dual\cite{Ma:rep}\cite{Ma:cat} or `centre' of
a monoidal category. The solution is the braided group
$\Aut(\CC)\rcocross
\Aut(\CC)^*$, which we show is dual-quasitriangular. This makes it an
example of some kind of braided codouble.

In Section~4 we make some remarks about double bosonisations and
general braided double crossproducts, also a topic of recent
interest. In fact, it was mentioned in the introduction of
q-alg/9511001\cite{Ma:dbos} that a theory of double cross products
$B\dcross C$ works fine in a braided category but does {\em not}
include the general construction of the `braided double' $B\dcross
B^*$ due to becoming `tangled up'.

\section{Reconstruction from $\CC^B\to \CC$}

Let $\CC$ be a braided category\cite{JoyStr:bra}, with braiding
$\Psi=\epsfbox{braid.eps}$ and inverse braiding
$\Psi^{-1}=\epsfbox{braidinv.eps}$. We use the diagrammatic methods
for braided groups due to the author in
\cite{Ma:bra}\cite{Ma:tra}\cite{Ma:bos}\cite{Ma:introm}, where the
product is represented as $\epsfbox{prodfrag.eps}$ etc. Contrary to
recent misconceptions, such a notation for braided algebra is {\em
not} in Yetter's fine paper\cite{Yet:rep} on crossed modules;
$\epsfbox{prodfrag.eps}$ there refers to ordinary Hopf algebras in
the category of vector spaces and $\epsfbox{braid.eps}$ to a
completely different braided category (so if one looks at the paper
in detail, there is nothing like braided groups or braided algebra
in \cite{Yet:rep}). We suppress the associativity $X\tens(Y\tens
Z)\isom (X\tens Y)\tens Z$ for objects $X,Y,Z\in
\CC$ by Mac Lane's coherence theorem. We assume throughout
that $\CC$ is rigid, i.e. every object comes with a dual $X^*$,
evaluation $\ev=\cup$ and coevaluation $\coev=\cap$. Here
$\cap:\und 1\to X\tens X^*$ but the unit object for the tensor
product is omitted in the diagrammatic notation. All categories are
assumed equivalent to small ones.

We let $\CC$ be such that the identity functor $i:\CC\to\CC$ obeys
the representability condition for comodules, i.e. there exists an
object $A\in \CC$ such that
\[ \theta:\Mor(A,V)\isom \Nat(i,i\tens V)\]
for any object $V$ in $\CC$ and such that the induced maps
\[ \theta_n:\Mor(A^n,V)\isom\Nat(i^n,i^n\tens V) \]
are also isomorphisms. Here the induced maps $\theta_n$ defined by
composing with the morphisms $\{\beta_X:X\to X\tens A\}$ making up
the natural transformation corresponding to $\id:A\to A$. In this
situation one has \cite{Ma:bg} a braided group $A=\Aut(\CC)$ living
in $\CC$ and $\beta_X$ is a {\em tautological coaction} of it on
each object $X$. For example, we may suppose for convenience that
$\CC$ is rigid and cocomplete over itself. These same constructions
go through more generally when, for example, $i:\CC\to\bar{\CC}$
where $\bar{\CC}$ is a larger category (eg some cocompletion of
$\CC$), yielding $\Aut(\CC)\in\bar{\CC}$). We use for that the more
general reconstruction from a functor $F:\CC\to\CV$, which is
defined similarly to the identity case above\footnote{And does
{\em not} require either $\CC$ or $\CV$ to be both cocomplete and rigid, 
see p205-206 of \cite{Ma:bg}. One may take for example $\CV$ cocomplete 
and the image of $\CC$ rigid, as explained in \cite{Ma:introm}, 
again long before \cite{Par:rec}.}; 
see \cite{Ma:bg}\cite{Ma:introm}\cite{Ma:book}. 

Now let $B$ be another braided group in $\CC$. The category $\CC^B$
of all braided comodules $(X,\cora)$ is monoidal\cite{Ma:bg}, where
$X\in\CC$ and $\cora:X\to X\tens B$ is the coaction.

\begin{lemma} In the above situation, the forgetful functor
$F:\CC^B\to \CC$ satisfies the representability assumption with
respect to the object $\Aut(\CC)\tens B$.
\end{lemma}
\proof We write $A=\Aut(\CC)$. Given any $\phi:A\tens B\to V$ we define
$\Theta(\phi)\in \Nat(F,F\tens V)$ by
$\Theta(\phi)_{(X,\cora)}=(\id\tens\phi)(\beta_X\tens\id)\cora$. It
is clearly a natural transformation since it can also be written as
$\Theta(\theta^B)_X=\theta^B_X\circ\cora$, where $\theta^B\in
\Nat(i\tens B,i\tens V)$ corresponds to $\phi$ via
$\theta^B(\phi)_X=(\id\tens\phi)\beta_X$.  In the converse
direction, given such a natural transformation $\Theta$, we define
$\theta^B(\Theta)_X=(\id\tens\eps\tens\id)\circ\Theta_{(X\tens
B,\Delta)}$ as clearly a natural transformation in $\Nat(i\tens
B,i\tens V)$, where we view $X\tens B$ as in $\CC^B$ by the trivial
coaction on $X$ and the coproduct on $B$. This then corresponds to
a morphism $\phi:A\tens B\to V$. These constructions are mutually
inverse. Thus, given $\Theta$ we define $\theta^B$ as stated. Then
\[\Theta(\theta^B)_X=\theta^B_X\circ\cora=(\id\tens\eps\tens\id)
\Theta_{(X\tens B,\Delta)}\circ\cora=(\id\tens\eps\tens\id)
(\cora\tens\id)\Theta_X=\Theta_X\] since $\cora:X\to X\tens B$ is a
morphism $(X,\cora)\to (X\tens B,\Delta)$ in $\CC^B$ (due to
$\cora$ a coaction) and $\Theta$ is natural. Conversely, given
$\theta^B$ we define $\Theta$ as stated. Then
\[ \theta^B(\Theta)_X=(\id\tens\eps\tens\id)
\Theta_{(X\tens B,\Delta)}=
(\id\tens\eps\tens\id)\theta^B_{X\tens
B}(\id\tens\Delta)=\theta^B_X\] since $\eps:X\tens B\to X$ is a
morphism in $\CC$ and $\theta^B$ is natural. Similarly for the
higher order representability conditions.
\endproof

\begin{theorem}cf. \cite{Ma:bos} Let $B$ be a braided group in $\CC$.
Braided reconstruction\cite{Ma:bg} from the forgetful functor
$F:\CC^B\to
\CC$ yields the {\em prebosonisation} braided group
$\Aut(\CC)\rcocross B$, with the braided tensor product algebra and
the cross coproduct coalgebra by the tautological coaction
$\beta_B:B\to B\tens \Aut(\CC)$ as an object of $\CC$.
\end{theorem}
\begin{figure}
\[\epsfbox{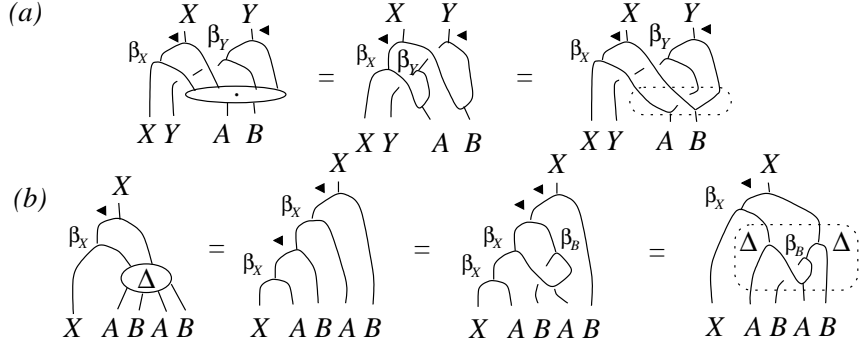}\]
\caption{Proof of Theorem~2.2}
\end{figure}\proof We routinely apply the reconstruction
theorem as presented diagrammatically in
\cite{Ma:introm}\cite{Ma:book}. The representing object is $A\tens
B$ from Lemma~2.1 and
$\beta_{(X,\cora)}=(\beta_X\tens\id)\circ\cora$ corresponds to the
identity on $A\tens B$. The product is defined in terms of this and
$\beta_{(X\tens Y,\cora\tens\cora)}$, which is the middle box in
Figure~1(a). From this we see that the product on $A\tens B$ is the
braided tensor product algebra (the dotted box). The coproduct is
defined as such that $\beta_{(X,\cora)}$ is a coaction, see the
first equality in Figure~1(b). The second equality is naturality
under $\cora:X\to X\tens B$ as a morphism in $\CC$. We then use
that $\cora$ is a coaction. This identifies the reconstructed
coproduct as a cross coproduct by $\beta_B$ (the dotted box). Cross
product braided groups are in \cite{Ma:bos} and we turn that
up-side-down. The result is a braided group due to the
braided-commutativity of $\Aut(\CC)$\cite{Ma:bg}.
\endproof

An example $BGL_q(2)\rcocross A_q^2$, where $B=A_q^2$ is the
quantum plane, is computed explicitly in \cite{Ma:com}. Note that
to obtain $BGL_q(2)$ (the braided group version of $GL_q(2)$) one
uses the braided reconstruction in the slightly more general form
where $\CC$ consists of finite-dimensional objects and $\Aut(\CC)$
lives more precisely in its cocompletion\cite{Ma:bg}.

Also, we can clearly iterate Theorem~2.2 to obtain braided groups
\[ (\Aut(\CC)\rcocross(\Aut(\CC)\rcocross\cdots\rcocross B)\cdots )\]
or `braided chain lattices'. Although there is in general no tensor
product of braided groups, Theorem~2.2 says that we can always
`tensor' by $\Aut(\CC)$ in this way.

Now we consider $B=\Aut(\CC)$ itself. By Theorem~2.2 we obtain a
braided group $\Aut(\CC)\rcocross\Aut(\CC)$, which we call the {\em
square} of $\Aut(\CC)$.

\begin{lemma} Let $B$ be any braided group. Then $B\rcocross B$
by the braided adjoint coaction (and braided tensor algebra) is
isomorphic to $B\lcross B$ by the braided adjoint action (and
braided tensor coalgebra).
\end{lemma}
\begin{figure}
\[ \epsfbox{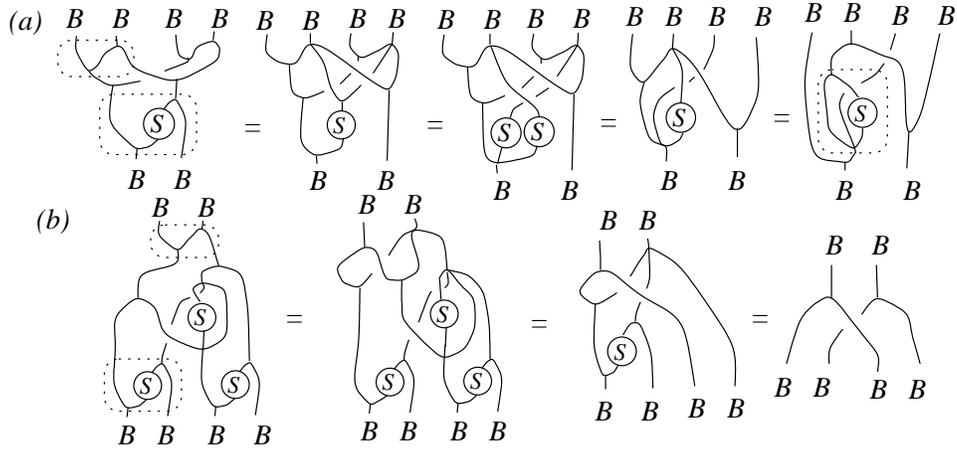}\]
\caption{Proof of Lemma~2.3}
\end{figure}
\proof In general, neither of these will
be braided groups -- so this is an isomorphism of algebras and
coalgebras. The former is shown in Figure~2(a). We apply the
isomorphism (the upper dotted box), then the braided tensor product
algebra, then the inverse of the isomorphism (lower dotted box). We
use the coproduct homomorphism property, antipode
antimultiplicativity, the antipode axiom to cancel a loop, and
finally recognise the cross product by the adjoint action (dotted
box on the right). Figure~2(b) does the computation for the
coproduct; we apply the isomorphism (dotted box), the cross
coproduct by the braided adjoint action, and then the inverse of
the isomorphism (dotted box). We use the coproduct homomorphism
property and cancel two resulting antipode loops. Cancelling the
resulting antipode loop gives the braided tensor coproduct on the
right. \endproof

\begin{propos} Reconstruction from the forgetful functor
$\CC^{\Aut(\CC)}\to\CC$ yields
\[ \Aut(\CC)\rcocross\Aut(\CC)\isom \Aut(\CC)\lcross\Aut(\CC)\]
where the left hand side is a cross coproduct by the adjoint
coaction (and braided tensor product algebra) and the left hand
side is the cross product\cite{Ma:bos} by the adjoint
action\cite{Ma:lie} (and braided tensor product coalgebra).
\end{propos}
\proof From the explicit realisation as a coend\cite{Ma:bg}
\cite{Lyu:tan}
$\Aut(\CC)=\int_X X^*\tens X$ it is clear that the canonical
coaction $\beta_{\Aut(\CC)}$ is the braided adjoint coaction
\cite{Ma:introm} of any
braided group $B$ on itself. Indeed, the coaction
$\id\tens\beta_X$ on each $X$ coincides with the coproduct
$\Aut(\CC)\to\Aut(\CC)\tens\Aut(\CC)$ as part of the
reconstruction,  the coaction on $X^*$ is conjugate to this via the
antipode $S$. We then use Lemma~2.3. \endproof

One knows from \cite{Ma:diag} that any trivial principal bundle has
a cross product form, and Proposition~2.4 tell us that in the
present case it is $\Aut(\CC)\lcross\Aut(\CC)$. This is then
explicitly a trivial braided principal bundle with right coaction
$\id\tens\Delta$ and the canonical inclusion of the left hand
$\Aut(\CC)$ as `base' of the bundle. It also means that the
Hopf algebra is as close as one can come to
$\Aut(\CC)\und\tens\Aut(\CC)$ as a braided group.

A concrete example is $BG_q\rcocross BG_q$ where the braided
coordinate rings $BG_q$ are quotients of the braided
matrices\cite{Ma:exa} $B(R)$. Writing the matrix generators of the
first copy as $\vecu$ and the second as $\vecv$, their relations
and that of the braided tensor product (braid statistics) between
them are
\[ R_{21}\vecu_1R\vecu_2= \vecu_2 R_{21} \vecu_1 R,\quad
R_{21}\vecv_1R\vecv_2= \vecv_2 R_{21} \vecv_1 R,\quad
R^{-1}\vecv_1R\vecu_2= \vecu_2 R^{-1}\vecv_1 R.\] The third
relations here correspond to the braiding\cite{Ma:exa}
\[ \Psi(R^{-1}\vecu_1\tens R\vecu_2)=\vecu_2 R^{-1}\tens \vecu_1 R\]
(written in \cite{Ma:exa} with all $R$ to one side) between any two
independent copies of $B(R)$. By Theorem~2.2, this $BG_q\rcocross
BG_q$ is a braided group with
\[ \Delta \vecu=\vecu\tens\vecu,\quad
\Delta \vecv=\vecu^{-1}\bullet \vecv\tens\vecu\vecv\]
where $\vecu^{-1}$ is to be exchanged with $\vecv$ using the
relations $R^{-1}\vecu_1\bullet R\vecv_2=
\vecv_2 R^{-1}\bullet\vecu_1 R$ and multiplied with $\vecu$ (these
are the relations of the braided tensor product of the copy
generated by $\vecv$ with that generated by $\vecu$, as a way of
describing the braided adjoint coaction as
conjugation\cite{Ma:lin}).

By Proposition~2.4, this is isomorphic as a braided group to
$BG_q\lcross BG_q$, with relations
\[ R_{21}\vecu_1R\vecu_2= \vecu_2 R_{21} \vecu_1 R,\quad
R_{21}\vecv_1R\vecv_2= \vecv_2 R_{21} \vecv_1 R,\quad
R_{21}\vecv_1R\vecu_2= \vecu_2 R_{21}\vecv_1 R\] and the coproduct
$\Delta\vecu=\vecu\tens\vecu$, $\Delta\vecv=\vecv\tens\vecv$. Also,
since $B(R)$ is also the braided enveloping bialgebra $U(\CL)$,
where $\CL$ is the braided Lie algebra associated to the
$R$-matrix\cite{Ma:lie}, we have a braided group $BG_q\lcross
U(\CL)$, which can be viewed as the algebra of observables of a
braided particle with `generalised momentum' $\CL$ moving under the
adjoint action on the braided space $BG_q$. The braided-Lie bracket
is the adjoint action and we use the known R-matrix formulae for
that to obtain
\[ \vecv_1R\vecu_2=\cdot\circ\vecv_1\la\Psi(\vecv_1\tens R\vecu_2)
=[\vecv_1,R\vecu_2]
R^{-1}\vecv_1R=R_{21}^{-1}\vecu_2R_{21}\vecv_1R\] to derive the
cross product as stated above. This braided group $BG_q\lcross
U(\CL)$ is related by transmutation to the quantum double as
explained in \cite{Ma:skl}. See also the next section.

These formulae also work fine at the braided bialgebra level
$B(R)\lcross B(R)$ for any biinvertible $R$-matrix, as one may
verfiy directly. And in spite of its origin as a braided cross
product, the formulae are remarkably symmetric between the two
copies of $B(R)$, reflecting the role as `tensor product'.
Moreover, the construction can be iterated to $n$ copies of $BG_q$
or $B(R)$, i.e. a `braided spin chain' $B(R)\lcross
B(R)\lcross\cdots\lcross B(R)$. This is generated by $\vecu^{(i)}$
(one for each copy of $B(R)$) with relations
\[ R_{21}\vecu_1^{(i)}R\vecu_2^{(j)}= \vecu_2^{(j)} R_{21}
\vecu_1^{(i)} R,\quad\forall i\ge j\]
and coproduct $\Delta \vecu^{(i)}=\vecu^{(i)}\tens
\vecu^{(i)}$, forming a braided group with braiding $\Psi$ as above
between any two copies of $B(R)$. This braided group is related by
transmutation to the iterated double crossproducts
$A(R)\dcross\cdots\dcross A(R)$ in \cite{Ma:mor}. It would be very
interesting to relate this approach also to the multiloop braided
algebras in \cite{Nil:str}.

\section{Reconstruction from $\CC^\circ$}

Given a monoidal category $\CC$ one has a dual monoidal category
$\CC^\circ$\cite{Ma:rep}\cite{Ma:cat}. It also arose at about the
same time as a `centre' or `double' construction, see
\cite{Ma:book}. Objects of $\CC^\circ$ are pairs $(V,\lambda)$
where $V$ is an object of $\CC$ and $\lambda_X:V\tens X\to X\tens
V$ is a natural transformation $\lambda\in\Nat(V\tens i,i\tens V)$
which `represents' the tensor product of $\CC$ in the sense
\[\lambda_{\und 1}=\id,\quad \lambda_{X\tens Y}
=\lambda_Y\circ\lambda_X.\]
Morphisms are morphisms of $\CC$ intertwining the corresponding
$\lambda$. Actually (an observation due to Drinfeld) $\CC^\circ$ in
the present case is braided, with
$\Psi_{(V,\lambda),(W,\mu)}=\lambda_W$. In our present applications
the category $\CC$ is itself braided as well.

\begin{propos} Reconstruction from the forgetful functor 
$\CC^\circ\to\CC$ yields the dual-quasitriangular braided group
$\Aut(\CC)\rcocross\Aut(\CC)^*$, a cross coproduct by the braided
coadjoint coaction.
\end{propos}
\begin{figure}
\[ \epsfbox{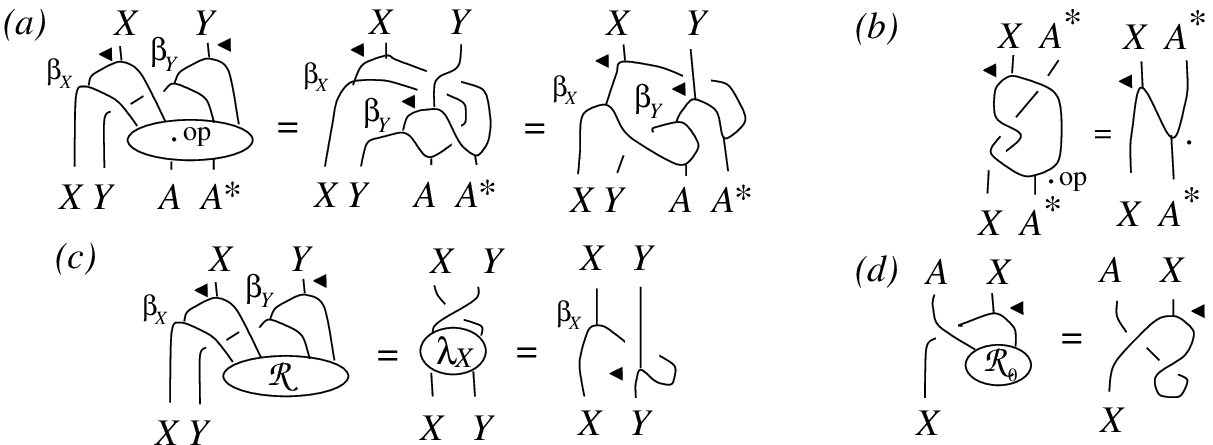}\]
\caption{Proof of Proposition~3.1}
\end{figure}
\proof From \cite{Ma:cat} we know that $\CC^\circ\isom \CC_{\Aut(\CC)}$,
the category of right $\Aut(\CC)$-modules, which is essentially the same
thing
as $\Aut(\CC)^*$-comodules (we assume that a suitable dual braided group
exists).
Then we can apply Theorem~2.2 with $B=\Aut(\CC)^*$ to obtain $\Aut(\CC)
\rcocross\Aut(\CC)^*$ as the reconstructed braided group. There is also a
second `opposite' product defined in Figure~3(a) via
$\beta_{(Y\tens X,\cora
\tens\cora)}$. We then use that $A=\Aut(\CC)$ itself is
braided-commutative.
If $A^*$ has an opposite product itself characterised by Figure~3(b) for
all $X$, then
we see that $(A\rcocross A^*)^{\rm op}=A\und\tens A^{*\rm op}$ (the braided
tensor product). Finally in the braided reconstruction theory there is a
dual-quasitriangular structure defined in Figure~3(c) by the braiding of
$\CC^\circ$, which we write with $\lambda$ in terms of the corresponding
coaction $\cora$ according to \cite{Ma:cat}. We see that if there is
morphism $\CR_0:A\tens A^*\to\und 1$ obeying Figure~3(d) then
$\CR=\CR_0\circ(\id\tens\eps\tens\eps
\tens\id)$. \endproof

This braided group $\Aut(\CC)\rcocross\Aut(\CC)^*$ is therefore
some kind of `braided codouble' of $\Aut(\CC)$ in spite of the the
fact that in general in a braided category this does not exist (see
the next section). Explicit examples are similar to those in the
preceding section since $BG_q$ is essentially self-dual.

\section{Braided double cross products}

This section is a kind of appendix to \cite{Ma:dbos}. It was
mentioned in its introduction:
\medskip

``the double cross product $B\dcross C$
construction does go through in a braided category, but the key
example of a general braided double $B\dcross B^{*\rm op}$ does
not''
\medskip

\noindent (a paraphrase) -- but details of the braided double cross
product were left unpublished due to this basic lack of
examples. Since that work, there has nevertheless been a lot of
interest in braided versions and generalisations of double cross
products and
bicrossproducts\cite{BesDra:cro}\cite{ZhaChe:dou}, and for this
reason we would like to publish now our calculations\cite{Ma:dcro}
mentioned in \cite{Ma:dbos}. They can by now be viewed as a special
case of the general constructions in
\cite{BesDra:cro}\cite{ZhaChe:dou}, but a case important enough to
study directly as we do now.

\begin{propos}\cite{Ma:dcro} A {\em matched pair} of braided groups
is $(B,C,\la,\ra)$ where  $B,C$ are braided groups in a category
$\CC$, $\la$ makes $B$ a braided left $C$-module
coalgebra\cite{Ma:introm}, $\ra$ makes $C$ a braided right
$B$-module coalgebra and $\la,\ra$ obey the conditions in
Figure~4(a1)-(a3) (and are trivial acting on 1). In this case there
is a {\em double cross product} braided group $B\dcross C$ with
product the dotted box in Figure~4(b) and the tensor product
coproduct.
\end{propos}
\begin{figure}
\[ \epsfbox{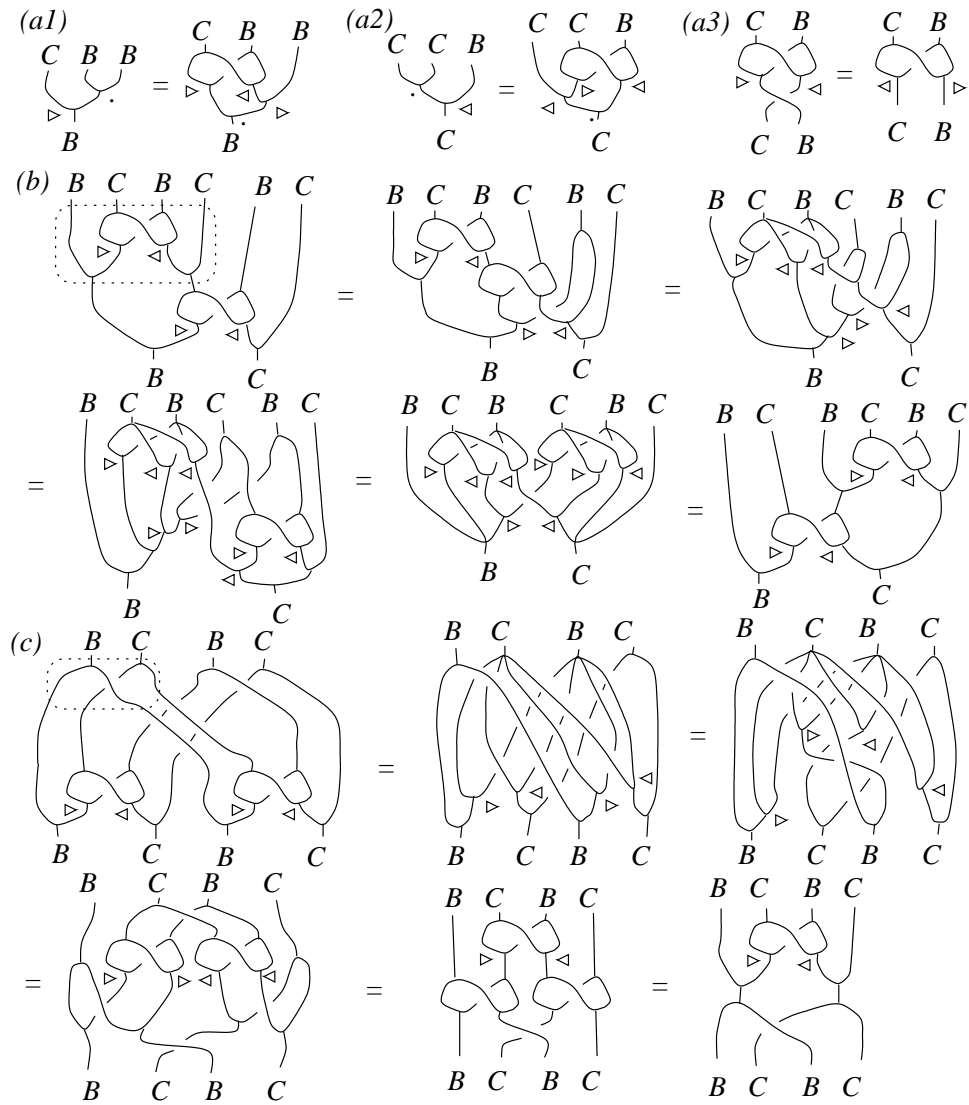}\]
\caption{Proof of Proposition~4.1}
\end{figure}
\proof Assuming a matched pair, Figure~4(b) checks that the braided
double cross product (the dotted box) is associative. We use the
coproduct homomorphism property for the first equality. The second
equality is that $\ra$ respects the coproduct of $C$ and that $\la$
is a left action. The third equality is axiom (a2). The fourth is a
reorganisation and associativity of $B,C$ to obtain an expression
which is symmetric under mirror-reflection (followed by reversal of
braid crossings). Therefore by the mirror image of the above steps,
taken in the reverse order, we obtain the final expression. Part
(c) checks that the braided tensor coproduct (the dotted box) obeys
the homomorphism property. The second equality is axiom (a3). The
third equality is another reorganisation. The fourth then uses that
$\la,\ra$ respect coproducts. Finally we use the coproduct
homomorphism properties of $B,C$.  \endproof

Conversely, given a braided group $X$ factorising into braided
groups $B,C$, one may recover $\la,\ra$ and $X\isom B\dcross C$ by
a similar proof to the unbraided case in \cite{Ma:book}.

As a possible example, we might try to build $D(B)=B\dcross B^{*\rm op}$
by braided coadjoint actions. The required coadjoint actions indeed
exist but axiom (a3) in Figure~4 fails due to `tangling up'. There
is no such problem in a {\em symmetric} monoidal category, however.
In this case, when $\CC={}_H\CM$ is the modules over a triangular
Hopf algebra, one obtains the bosonisation of this
symmetric-category double as
\[ (B\dcross B^{*\rm op})\lbiprod H\isom B\lbiprod
H\rbiprod B^{*\rm op},\]
the double-bosonisation. This is explained in \cite{Ma:pra96} (in
the super case) and was indeed one of the main motivations behind
the double-bosonisation construction in \cite{Ma:dbos}; in general
$B\dcross B^{*\rm op}$ does not exist but the double-bosonisation
does!


\begin{thebibliography}{10}

\bibitem{Ma:bg}
S.~Majid.
\newblock Braided groups.
\newblock {\em J. Pure and Applied Algebra}, 86:187--221, 1993.

\bibitem{Ma:bra}
S.~Majid.
\newblock Braided groups and algebraic quantum field theories.
\newblock {\em Lett. Math. Phys.}, 22:167--176, 1991.

\bibitem{Ma:tra}
S.~Majid.
\newblock Transmutation theory and rank for quantum braided groups.
\newblock {\em Math. Proc. Camb. Phil. Soc.}, 113:45--70, 1993.

\bibitem{Ma:bos}
S.~Majid.
\newblock Cross products by braided groups and bosonization.
\newblock {\em J. Algebra}, 163:165--190, 1994.

\bibitem{Ma:rep}
S.~Majid.
\newblock Representations, duals and quantum doubles of monoidal categories.
\newblock {\em Suppl. Rend. Circ. Mat. Palermo, Ser. II}, 26:197--206, 1991.


\bibitem{Ma:cat}
S.~Majid.
\newblock Braided groups and duals of monoidal categories.
\newblock {\em Canad. Math. Soc. Conf. Proc.}, 13:329--343, 1992.

\bibitem{Ma:exa}
S.~Majid.
\newblock Examples of braided groups and braided matrices.
\newblock {\em J. Math. Phys.}, 32:3246--3253, 1991.

\bibitem{Ma:lin}
S.~Majid.
\newblock Quantum and braided linear algebra.
\newblock {\em J. Math. Phys.}, 34:1176--1196, 1993.

\bibitem{Ma:skl}
S.~Majid.
\newblock Braided matrix structure of the Sklyanin algebra and of
the quantum Lorentz group.
\newblock {\em Commun. Math. Phys.} 156:607--638, 1993.

\bibitem{Ma:lie}
S.~Majid.
\newblock Quantum and braided {L}ie algebras.
\newblock {\em J. Geom. Phys.}, 13:307--356, 1994.

\bibitem{Ma:diag}
S.~Majid.
\newblock Diagrammatics of braided group gauge theory.
\newblock {\em Preprint}, 1996.

\bibitem{Ma:dbos}
S.~Majid.
\newblock Double bosonisation and the construction of {$U_q(g)$}, 1995.
\newblock To appear in {\em Math. Proc. Camb. Phil. Soc.}

\bibitem{Ma:introm}
S.~Majid.
\newblock Algebras and {H}opf algebras in braided categories.
\newblock volume 158 of {\em Lec. Notes in Pure and Appl. Math},
pages 55--105.
  Marcel Dekker, 1994.

\bibitem{Ma:book}
S.~Majid.
\newblock {\em Foundations of Quantum Group Theory}.
\newblock Cambridge Univeristy Press, 1995.

\bibitem{Par:rec}
B.~Pareigis.
\newblock Reconstruction of hidden symmetries.
\newblock {\em J. Algebra}, 183:90--154, 1996.

\bibitem{JoyStr:bra}
A.~Joyal and R.~Street.
\newblock Braided monoidal categories.
\newblock Macquarie Report, 1986.

\bibitem{Yet:rep}
D.N. Yetter.
\newblock Quantum groups and representations of monoidal categories.
\newblock {\em Math. Proc. Camb. Phil. Soc.}, 108:261--290, 1990.

\bibitem{Ma:com}
S.~Majid.
\newblock Some comments on bosonisation and biproducts.
\newblock {\em Czech J. Phys.}, 47:151--171, 1997.


\bibitem{Lyu:tan}
V.V. Lyubashenko.
\newblock Tangles and {H}opf algebras in braided categories.
\newblock {\em J. Pure and Applied Algebra}, 98:245--278, 1995.

\bibitem{Ma:mor}
S.~Majid.
\newblock More examples of bicrossproduct and double cross product
Hopf algebras.
\newblock {\em Isr. J. Math}, 72:133--148, 1990.

\bibitem{Nil:str}
\newblock F.~Nill.
\newblock Structure of monodromy algebras and Drinfeld doubles.
\newblock {\em Rev.Math.Phys.}, 9:371-395, 1997.

\bibitem{BesDra:cro}
Yu.~N. Bespalov and B.~Drabant.
\newblock Cross product bialgebras, I.
\newblock {\em Preprint}, Damtp/98-9, February, 1998.

\bibitem{ZhaChe:dou}
S.~Zhang and H.-X. Chen.
\newblock The double biproduct in braided tensor categories.
\newblock {\em Preprint}, Fudan/Nanjing Univ., 1998.

\bibitem{Ma:dcro}
S.~Majid.
\newblock Notes to ref. [12], {\em Unpublished}, 1995.

\bibitem{Ma:pra96}
S.~Majid.
\newblock New quantum groups by double-bosonisation.
\newblock {\em Czech J. Phys.}, 47:79--90, 1997.

\end{thebibliography}

\end{document}